\documentclass[a4paper,11pt]{amsart}

\usepackage[utf8]{inputenc}
\usepackage[english]{babel}
\usepackage{amsmath}
\usepackage{amsthm}
\usepackage{amsfonts}
\usepackage{amssymb}
\usepackage[autostyle]{csquotes} 
\usepackage{caption}
\usepackage{hyperref} 
\usepackage[english]{cleveref} 
\usepackage{diagbox}
\usepackage{dsfont}
\usepackage{emptypage}
\usepackage[shortlabels]{enumitem}
\usepackage{extpfeil}
\usepackage{etoolbox}
\usepackage{faktor}
\usepackage{fancyhdr}
\usepackage{float}
\usepackage{graphicx}
\usepackage{indentfirst}
\usepackage{mathdots}
\usepackage{mathtools}
\usepackage{mathrsfs}
\usepackage{multirow}
\usepackage{rotating}
\usepackage{stackrel}
\usepackage{stmaryrd}
\usepackage{subcaption}
\usepackage{tikz}
\usepackage{verbatim}

\usetikzlibrary{patterns,arrows,calc,matrix,decorations.pathreplacing}

\newcommand*\T{\mathcal{T}}
\newcommand*\G{\mathcal{G}}
\newcommand*\E{\mathcal{E}}
\newcommand*\R{\mathbb{R}}
\newcommand*\Z{\mathbb{Z}}

\newcommand*\Q{\mathbb{Q}}
\newcommand*\C{\mathbb{C}}

\newcommand{\defeq}{\stackrel{\textnormal{def}}{=}}

\theoremstyle{definition}

\theoremstyle{remark}

\DeclareMathOperator{\dd}{d}
\DeclareMathOperator{\MM}{M}
\DeclareMathOperator{\M}{\mathcal{M}}
\DeclareMathOperator{\rk}{rk}
\newcommand{\ie}{i.e.\ }

\makeatletter
\newcommand*{\bigcdot}{%
  {\mathbin{\mathpalette\bigcdot@{}}}%
}
\newcommand*{\bigcdot@scalefactor}{.75}
\newcommand*{\bigcdot@widthfactor}{1.4}
\newcommand*{\bigcdot@}[2]{%
  \sbox0{$#1\vcenter{}$}
  \sbox2{$#1\cdot\m@th$}%
  \hbox to \bigcdot@widthfactor\wd2{%
    \hfil
    \raise\ht0\hbox{%
      \scalebox{\bigcdot@scalefactor}{%
        \lower\ht0\hbox{$#1\bullet\m@th$}%
      }%
    }%
    \hfil
  }%
}
\makeatother

\robustify{\bigcdot}

\begin{document}
\title[Poincaré polynomial of elliptic arrangements]{Poincaré polynomial of elliptic arrangements is not a specialization of the Tutte polynomial}
\author[R. Pagaria]{Roberto Pagaria}
\address{Roberto Pagaria}
\address{\textup{Scuola Normale Superiore\\ Piazza dei Cavalieri 7, 56126 Pisa\\ Italia}}
\email{roberto.pagaria@sns.it}
\begin{abstract}
The Poincaré polynomial of the complement of an arrangements in a non compact group is a specialization of the $G$-Tutte polynomial associated with the arrangement.
In this article we show two unimodular elliptic arrangements (built up from two graphs) with the same Tutte polynomial, having different Betti numbers.
\end{abstract}

\maketitle

\section{Introduction}

Let $A \in \MM(k,n;\Z)$ be an integer matrix and let $G$ be a group of the form $H \times (S^1)^p \times \R^q$, where $H$ is a finite abelian group.
Each column $\alpha$ of $A$ defines a morphism from $G^k$ to $G$ given by
\[(g_1, \dots, g_k) \mapsto \alpha_1 g_1 + \alpha_2 g_2 + \dots + \alpha_k g_k .\]
We call $H_i \subset G^k$ the kernel of the map defined by the $i^{\textnormal{th}}$-column of $A$.
The \textit{complement} of the \textit{arrangement} $A$ in $G$ is the topological space
\[ \M(A;G)= G^k \setminus \bigcup_{i=1}^n H_i.\]

When $G=\R^2\simeq \C$ we obtain the classical definition of \textit{hyperplane arrangements}.
If $G=S^1 \times \R \simeq \C^*$ the arrangement is called \textit{toric}.
We are mainly interested in the case $G=S^1 \times S^1 \simeq E$ (an elliptic curve), this arrangements is called \textit{elliptic arrangement}.

There are several combinatorial objects associated with an arrangement: for instance, the poset of layers, the arithmetic matroid (\cite{BrandenMoci2014,DM}) and the $G$-Tutte polynomial (\cite{MociTutte2012, LTY2017, TY2018}).
Given a subset $S$ of $[n]=\{1, 2, \dots, n\}$ we call \textit{layer} any connected component of the intersection $\bigcap_{i \in S} H_i$.
The \textit{poset of layers} is the set of all layers ordered by reverse inclusion.
The \textit{arithmetic matroid} is the triple $([n], \rk, m_G)$ associated with toric, hyperplanes or elliptic arrangements, where $\rk(S)$ and $m_G(S)$ are, respectively, the codimension and the number of connected components of $\bigcap_{i \in S} H_i$.
The $G$\textit{-Tutte polynomial} is a generalization of the arithmetic Tutte polynomial and of the classical Tutte polynomial; it is defined by
\[ T_A^G(x,y) \defeq \sum_{S\subseteq [n]} m_G(S) (x-1)^{\rk [n] - \rk(S)}(y-1)^{|S|-\rk(S)}.\]

Recently, a formula for the Poincaré polynomial of $\M(A;G)$ was found by Liu, Tran and Yoshinaga \cite{LTY2017} when $G$ is not compact, \ie $q>0$.
This formula involves the $G$\textit{-characteristic polynomial} $\chi_A^G(t)$, which is a specialization of the $G$-Tutte polynomial:
\[ \chi_A^G(t)= (-1)^{\rk [n]} t^{k-\rk [n]} T_A^G(1-t,0).\]
When $G$ is not compact, the Poincaré polynomial of $\M(A;G)$ is
\[ P_{\M(A;G)}=(-t^{p+q-1})^k \chi_A^G \left(- \frac{P_G(t)}{t^{p+q-1}}\right),\]
where $P_G(t)=m_G(\emptyset)(t+1)^p$ is the Poincaré polynomial of the group $G$.
The formula 
\[e(\M(A;G)) = (-1)^{(p+q)k}\chi_A^G((-1)^{p+q} e(G)) \]
for the Euler characteristic holds for all groups $G$ ($e(G)$ is the euler characteristic of $G$), see \cite{Bibby2016,LTY2017}.

We focus on the ``smallest'' compact group $G= S^1 \times S^1$, the case $G=S^1$ being trivial.
From now on, we denote the two-dimensional compact torus $S^1 \times S^1$ by $E$.
In this case, Bibby \cite{Bibby2016} and Dupont \cite{DupontHypersurfaces2015} have given a model of the cohomology ring $H^\bullet(M(A;E);\Q)$, provided by the second page of the Serre spectral sequence for the inclusion $M(A;E) \hookrightarrow E^k$.
As shown in \cite{PagariaEC}, this model is combinatorial, \ie can be defined from the arithmetic matroid $([n],\rk,m_E)$.
Thus the Betti numbers are implicitly codified in the arithmetic matroid, but there is no explicit formula that allows their calculation.
We will show that these Betti number are independent from the arithmetic Tutte polynomial, exhibiting an example.

\section{The model for cohomology}
We recall the model developed by Dupont \cite{DupontHypersurfaces2015} and Bibby \cite{Bibby2016} for the cohomology ring in the particular case of graphic elliptic arrangements.

Let $E^{k+1}/E\simeq E^k$ be the quotient of $E^{k+1}$ by the diagonal action of $E$.
Given a finite graph $\G=([k+1],\E)$, undirected and without loops or multiple edges, we can define an arrangement $A_\G$ in $E^{k+1}/E$ given by the divisors
\[ H_e=H_{i,j} \defeq \{\underline{g}\in E^{k+1}/E \mid g_i=g_j\},\]
for each edge $e=(i,j)\in \E$.
We fix arbitrarily a spanning forest $\T$ of $\G$ and an orientation of $\G$.

Consider the external algebra $\Lambda$ over the rationals on 
the generators $\{\omega_{e}, x_{t},y_{t} \}_{\substack{e \in \E \\t \in \T}}$.
We set the bi-degree of  $\omega_{e}$ to be $(1,0)$ and the one of $x_{t}$ and $y_t$ to be $(0,1)$.
For the sake of an easier notation we define for each oriented edge $e= i \rightarrow j$ the element $x_e= \sum_{f \in \gamma(e)} \epsilon(e,f) x_f$, where $\gamma(e)$ is the unique path from $i$ to $j$ in $\T$ and $\epsilon(e,f)$ is $1$ if the arc $f$ is oriented as in the path $\gamma(e)$, $-1$ otherwise.
Consider the homogeneous ideal $I\subset \Lambda$ generated by the elements $\omega_e x_e$, $\omega_e y_e$ and
\[ \sum_{i=0}^l (-1)^i \omega_{e_0} \omega_{e_1} \dots \widehat{\omega}_{e_i} \dots \omega_{e_l}  ,\]
for every cycle $C=(e_0, e_1, \dots ,e_l)$ of $\G$.
We call $E_2(A_\G)$ the quotient $\Lambda/I$.
Finally, we define the differential $\dd_2 \colon E_2(A_\G) \to  E_2(A_\G)$ on the generators by $\dd_2(\omega_e)=x_ey_e$ and $\dd_2(x_e)= \dd_2(y_e)=0$.
This is well defined since $\dd_2(I) \subseteq I$.
The model $(E_2(A_\G), \dd_2)$ coincides with the second page of the Leray spectral sequence and the cohomology of the second page (\ie the third one) is the cohomology ring of $\M(A_\G;E)$ with rational coefficients.
The bi-gradation of the third page corresponds to the bi-gradation given by the mixed Hodge structure (and the total degree).
Let $e(a,b)$ be the dimension of the homogeneous subspace of bi-degree $(a,b)$ of the third page.
The number $e(a,b)$ coincides with the dimension of the subspace of $H^{a+b}(M(A_\G;E))$ of weight $a+2b$ (see \cite[4 pg~81]{Deligne75}).

Since the elliptic arrangement $A_\G$ is \textit{unimodular}, \ie every subset of divisors has connected intersection, the $G$-Tutte polynomial $T_{A_\G}^E$ coincides with the classical Tutte polynomial $T_\G$ associated with the graph $\G$.
In particular the dimension of $E_2^{a,b}(A_\G)$
can be easily calculated from 
\[ \sum_{a,b} \dim E_2^{a,b}(A_\G) t^a s^b= T_\G \left(1+\frac{(1+t)^2}{s},0 \right)s^k.\]
Thus, the Hodge polynomial evaluated in $(-1,u)$ is 
\[\sum_n \left( \sum_{m\geq 0} (-1)^m e(n-2m,m) \right)u^n =  T_\G \left(1-\frac{(1+u)^2}{u^2},0 \right)(-u^2)^k,\]
and the Euler characteristic of $\M(A_\G;E)$ is $(-1)^k T_\G(1,0)$, as shown in \cite{Bibby2016}.

\section{The example}

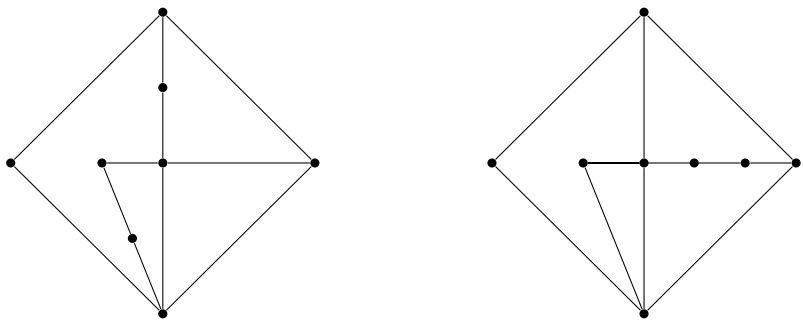
\begin{figure}
\centering
\begin{subfigure}{0.49\textwidth}
\centering
\begin{tikzpicture}[myvert/.style={fill,circle, minimum size = 0.12 cm, inner sep = 0 cm}]
\node[myvert] (g) at (0,2){};
\node[myvert] (h) at (2,0){};
\node[myvert] (a) at (-2,0){};
\node[myvert] (b) at (0,-2){};
\node[myvert] (e) at (0,0){};
\node[myvert] (d) at (-0.8,0){};
\node[myvert] (c) at (-0.4,-1){};
\node[myvert] (f) at (0,1){};
\draw[thin] (g) -- (a) -- (b) -- (c) -- (d) -- (e)-- (f) -- (g) -- (h) -- (b) -- (e) -- (h);
\end{tikzpicture}
\end{subfigure}
\begin{subfigure}{0.49\textwidth}
\centering
\begin{tikzpicture}[myvert/.style={fill,circle, minimum size = 0.12 cm, inner sep = 0 cm}]
\node[myvert] (g) at (0,2){};
\node[myvert] (h) at (2,0){};
\node[myvert] (a) at (-2,0){};
\node[myvert] (b) at (0,-2){};
\node[myvert] (e) at (0,0){};
\node[myvert] (d) at (-0.8,0){};
\node[myvert] (c) at (0.66,0){};
\node[myvert] (f) at (1.33,0){};
\draw[thin] (g) -- (a) -- (b) -- (d) -- (e) -- (g) -- (h) -- (b) -- (e) -- (d) -- (f) -- (h);
\end{tikzpicture}
\end{subfigure}
\caption{The graph $\G_1$ on the left and $\G_2$ on the right}
\label{fig:i_grafi}
\end{figure}

Consider the two graphs $\G_1$ and $\G_2$ in \cref{fig:i_grafi} and the corresponding graphic elliptic arrangements $A_1$ and $A_2$.
These graphs appeared for the first time in~\cite{Schwarzler}.
They share the same Tutte polynomial, which is the following
\[
\begin{split}
T(x,y)=x^7 + 4x^6 + x^5y + 9x^5 + 6x^4y + 3x^3y^2 + x^2y^3 + 13x^4 + 13x^3y + \\
+7x^2y^2 + 3xy^3 + y^4 + 12x^3 + 15x^2y + 9xy^2 + \\
+ 3y^3 + 7x^2 + 9xy + 4y^2 + 2x + 2y.
\end{split}
\]
Using SAGE \cite{sagemath}, we have computed the mixed Hodge numbers of $\M(A_1)$ and of $\M(A_2)$ and reported them in \Cref{tab:E3_di_G1,tab:E3_di_G2}.
For this computation we have used the code available \href{http://poisson.phc.dm.unipi.it/~pagaria/Graphic_Elliptic_Arr.txt}{here}\footnote{\url{http://poisson.phc.dm.unipi.it/~pagaria/Graphic_Elliptic_Arr.txt}}; the calculation of the Hodge number $e(4,2)$ has taken more than $2$ days with a CPU of $2.2$GHz and about $32$ GB of RAM.
Some Hodge numbers have being calculated using the following formula
\begin{align*}& \sum_n \left( \sum_{m\geq 0} (-1)^m e(n-2m,m) \right)u^n =  1 + 14 u + 80 u^2 + 232 u^3 + 329 u^4 +\\
& +122 u^5 - 165 u^6 - 24 u^7 + 
 164 u^8 - 56 u^9 - 71 u^{10} + 68 u^{11} - 26 u^{12} + 4 u^{13}.
 \end{align*}
The Poincaré polynomials of $\M(A_1)$ and $\M(A_2)$ are different:
\begin{align*}
P_{\M(A_1;E)} (t) &= 1+14 t+ 82t^2+269t^3+570t^4+820 t^5 +765t^6 + 363t^7, \\
P_{\M(A_2;E)} (t) &= 1+14 t+ 82t^2+270t^3+578t^4+844 t^5 +785 t^6 + 366 t^7.
\end{align*}
The Euler characteristic of $M(A_1;E)$ and of $M(A_2;E)$ are both equal to $-48$.
\begin {table}
\caption {The Hodge numbers of $\M(A_1)$}\label{tab:E3_di_G1} 
\begin{tabular}{cccccccc}
0 & 4 &  &  &  &  &  &  \\ 
0 & 6 & 26 &  &  &  &  &  \\ 
0 & 4 & 45 & 74 &  &  &  &  \\ 
0 & 8 & 54 & 154 & 116 &  &  &  \\ 
0 & 6 & 60 & 200 & 259 & 94 &  &  \\ 
0 & 2 & 29 & 144 & 302 & 224 & 41 &  \\ 
1 & 14 & 80 & 234 & 358 & 260 & 77 & 8 \\ 
\end{tabular}
\medskip
\caption*{In position $(i,j)$ there is the dimension of the subgroup of weight $i+2j$ in $H^{i+j}(\M(A_1;E))$.}
\end {table}

\begin{table}
\caption {The Hodge numbers of $\M(A_2)$}\label{tab:E3_di_G2} 
\begin{tabular}{cccccccc}
0 & 4 &  &  &  &  &  &  \\ 
0 & 6 & 26 &  &  &  &  &  \\ 
0 & 6 & 45 & 74 &  &  &  &  \\ 
0 & 10 & 69 & 162 & 116 &  &  &  \\ 
0 & 6 & 59 & 202 & 271 & 100 &  &  \\ 
0 & 2 & 30 & 150 & 301 & 224 & 38 &  \\ 
1 & 14 & 80 & 234 & 359 & 266 & 77 & 8 \\ 
\end{tabular}
\medskip
\caption*{In position $(i,j)$ there is the dimension of the subgroup of weight $i+2j$ in $H^{i+j}(\M(A_2;E))$.}
\end{table}

\subsection*{Acknowledgements}
The server, used for computations, has been acquired thanks to the support of the University of Pisa, within the call ``Bando per il cofinanziamento dell’acquisto di medio/grandi attrezzature scientifiche 2016''.

\bibliography{biblio_Poin}
\bibliographystyle{plain}
\bigskip
\end{document}